\documentclass{article}
\usepackage{bbm}
\usepackage{amssymb}
\usepackage{amsmath}
\usepackage{pxfonts}
\begin{document}

\begin{Large}
\centerline{\bf On Larcher's theorem concerning good lattice points}
\centerline{\bf 
 and}
\centerline{\bf
multiplicative subgroups modulo p.}
 \vskip+1.0cm

\centerline{{ N.G. Moshchevitin, D.M. Ushanov}
\footnote{research is supported  by RFBR grant No 09-01-00371a}}
\end{Large}
\vskip+1.0cm

\centerline{\bf Abstract.}
\begin{small}
We prove the existence of two-dimensional good lattice points in thick multiplicative subgroups modulo $p$.

\end{small}
\vskip+0.7cm

\section{Introduction}

Let $p\ge 3$ be a prime number. Take an integer $a$ such that $1\le a \le p-1$. Consider a sequence of points 
\begin{equation}
\xi_x = \left( \frac{x}{p}, \left\{\frac{ax}{p}\right\} \right) \in [0,1]^2,\,\,\,\, x=0,1,2,...,p-1.
\label{point}
\end{equation}
Let 
$$N_p(\gamma_1,\gamma_2)=\#\{x : 0\leq x < p,\ \xi_x\in [0,\gamma_1]\times[0,\gamma_2]\}$$
and let 
$$D_p(a) = \sup_{\gamma_1,\ \gamma_2 \in [0,1]} |N_p(\gamma_1,\gamma_2) - \gamma_1\gamma_2 p|$$
be the discrepancy of the set (\ref{point}).

In \cite{1} G. Larcher proved a series of results on the existence of well-distributed sets of the form (\ref{point}).
For example he proved the existence of $a \in [0,1,\cdots,p-1]$
such that
$$
D_p(a) \leq c\log{p}\log\log{p}
$$
with an absolute constant $c.$

In the present paper we generalize this result.

In the sequel 
$\mathbb{Z}_p^*$ denotes the multiplicative group of residues modulus $p$.
$U$ denotes a multiplicative subgroup of $\mathbbm{Z}_p^*$ and $\|\cdot \|$ denotes the distance 
to the nearest integer.

For $1\leq a < p$ we need the continued fraction expansion
\begin{equation}
\frac{a}{p} = \cfrac{1}{b_1(a) + \cfrac{1}{ b_2(a) + \cdots + \cfrac{1}{b_l(a)} }},\ \ \ l = l(a).
\label{CF}
\end{equation}

\textbf{ Theorem 1. } \, {\it Let $p$ be prime, $U$ be a multiplicative  subgroup in $\mathbbm{Z}_p^*.$ For  $v \neq 0$ we consider the set $R = v\cdot U$ and let
$$
\#R \geq 10^5 p^{7/8}\log^{3/2}{p}.
$$
Then for at least a half of elements $a \in R$ all partial
quotients $b_j(a)$ in the continued fraction expansion (\ref{CF})
are less than $16\log{p}.$
}

Theorem 1 improves a result from \cite{3}.

\textbf{ Theorem 2. } \,   {\it Let $p$ be prime, $U$ be a multiplicative  subgroup in $\mathbbm{Z}_p^*.$  For  $v \neq 0$ we consider the set $R = v\cdot U$ and let
$$
\#R \geq 10^8 p^{7/8}\log^{5/2}{p}.
$$
Then there exists an element $a \in R,$ $a/p=[b_1, b_2, \cdots, b_l],$ $b_i = b_i(a),$ $l = l(a)$ with
$$
\sum_{i=1}^{l} b_i \leq 500\log{p}\log\log{p}.
$$
}

It is well known (see \cite{4}) that
$$
D_p(a) \ll
\sum_{1\leq i \leq l(a)}b_i(a).
$$
  So we immediately obtain the following

\textbf{ Corollary.}\,{\it
Under the conditions of Theorem 2 
  there exists an element $a \in R$ such that
$$
D_p(a) \ll \log{p}\log\log{p}.
$$
}

We 
do not calculate optimal constants in our results. Of course constants $10^5$ and $10^8$ may be reduced.
\section{Character sums}

Let $p$ be prime, $1 < t \leq p,$
$k = \sqrt{\frac{2p}{t},}$ $j = \lceil \log_2 \frac{p}{k} \rceil.$
Define rectangles
$$
\begin{array}{l}
\Pi_0 = [1,k]\times [1,k],\\
\Pi_1 = [k+1,2k]\times [1,k/2,]\\
\Pi_2 = [2k+1,4k]\times [1,k/4,]\\
\cdots,\\
\Pi_{\nu} = [2^{\nu-1}k + 1, 2^{\nu}k]\times [1,k/2^{\nu}],\\
\cdots;\\
\Pi_{-1} = [1,k/2]\times [k+1,2k],\\
\Pi_{-2} = [1,k/4]\times [2k+1,4k],\\
\cdots,\\
\Pi_{-\nu} = [1,k/2^{\nu}]\times [2^{\nu-1}k + 1, 2^{\nu}k],\\
\cdots,
\end{array}
$$
and let $\Pi^t = \cup_{i=-j}^{j} \Pi_{i},$ so $\Pi^t$ consists of $\leq 2 \log_2 p$ rectangles $\Pi_i.$ It's clear that
\begin{equation}
\{(x,y)\in \mathbb{Z}^2 \mid 1\leq x < p,\ 1\leq y<p,\ xy \leq p /t\}
\subset \Pi^t.\label{Pi}
\end{equation}
Moreover, for different $\nu$ and $\mu$ we have
$$
\Pi_\nu \cap \Pi_\mu = \varnothing.
$$

\textbf{ Lemma 1. } {\it Let $p$ be prime, $c \geq 1,$ $k = \sqrt{\frac{2p}{c},}$ $\chi$ be a non-principal character to prime modulus $p.$ Then
$$
\left| \sum_{ (x,u) \in \Pi^c } \chi(x) \overline{\chi(u)} \right| \leq
10000 p^\frac{7}{8} \log^2{p} / \sqrt{c}.
$$
}
\textit{ Proof. } Dividing the summation area into parts, we obtain
$$
\left| \sum_{ (x,u) \in \Pi^c } \chi(x) \overline{\chi(u)} \right| \leq
\sum_{i=-j}^j \left|\sum_{(x,u)\in \Pi_i} \chi(x) \overline{\chi(u)} \right|.
$$

Let $h$ denotes the height of rectangle $\Pi_i$ and $w$ denotes the width.
Then $hw = k^2.$

We will use following Burgess' result (see [4] for details)

\textbf{ Theorem. } {\it Let $\chi$ be a non-principal character to prime modulus. Then
$$
\left| \sum_{1 \leq x \leq N} \chi(x) \right| \leq 30N^{1 - \frac 1 r}p^{\frac{r+1}{4r^2}} (\log{p})^{\frac{1}{r}}.
$$
Here $r$ is an arbitrary positive integer. }

Taking $r=2$ in the Burgess' theorem we obtain
$$
\left|\sum_{(x,u)\in \Pi_i} \chi(x) \overline{\chi(u)} \right|
\leq 900 \sqrt{hw} p^\frac{3}{8} \log{p}
\leq 900 k p^\frac{3}{8} \log{p} = 900 \sqrt{\frac{2p}{c}}
p^\frac{3}{8} \log{p}.
$$

Since there is only $\leq 2\log_2{p}$ rectangles $\Pi_i$
we deduce that
$$
\left| \sum_{ (x,u) \in \Pi^c } \chi(x) \overline{\chi(u)} \right| \leq
10000 p^\frac{7}{8} \log^2{p} / \sqrt{c}
$$

and the lemma follows.

\section{Continued fractions}

We will use two lemmas about continued fractions (see \cite{5}).

\textbf{ Lemma A. } {\it If $\alpha \in \mathbbm{R},$ \ $\frac{a}{b} \in \mathbbm{Q}$ and
$$
\left| \alpha - \frac{a}{b}\right| < \frac {1}{2b^2},
$$
then $\frac{a}{b}$ is a convergent to $\alpha.$}

\textbf{ Lemma B. } {\it If $\frac{p_n}{q_n} \neq \alpha$ is the $n$-th convergent to $\alpha$ then
$$
\frac{1}{q_n(q_n+q_{n+1})} < \left| \alpha - \frac{p_n}{q_n} \right| < \frac{1}{q_n q_{n+1}}.
$$
}

\section{ Proof of Theorem 1}

Let $t = 16\log{p}.$ Consider the sum
$$
S(a) = \sum \delta_p(axy^* - 1),
$$
where
$$
\delta_p(x) =
 \begin{cases}
1, & \text{$ x\equiv 0 \pmod{p} $}\\
0, & \text{$ x\nequiv 0 \pmod{p} $}
\end{cases}
$$
$y^*\in \mathbb{Z}_p^*
$ is defined from
$$
yy^*\equiv 1\pmod{p},
$$
and the summation is over all pairs $(x,|y|) \in \Pi^t.$

If $S(a) = 0$ then by (\ref{Pi}) for all $1\leq x,$ $1\leq|y|<p,$ $x|y| \leq p/t$ we have
$$
ax-y \nequiv 0 \pmod{p}.
$$
Hence if
$$
\begin{cases}
ax-y\equiv 0 \pmod{p},\\
1\leq x < p,\\
1\leq |y| <p,
\end{cases}
$$
then
\begin{equation}
x|y|>\frac{p}{t}.
\label{st}
\end{equation}
In particular (\ref{st}) holds for $y = \pm \left\|\frac{ax}{p}\right\|p.$ Therefore
for each $1\leq x < p$ we have
\begin{equation}
\left\|\frac{ax}{p}\right\| > \frac{1}{xt}.
\label{a}
\end{equation}
Using Lemma B we obtain
\begin{equation}
\left\|\frac{ax_{k-1}}{p}\right\| < \frac{1}{x_{k-1}b_k(a)}.
\label{b}
\end{equation}
From (\ref{a}) and (\ref{b}) we see that $b_k(a) < t$ and all continued fraction coefficients of $a/p$ are less than
$t.$

Now we express $S(a)$ as a sum
$$
S(a) = \frac{1}{p-1}\sum_{\chi \pmod{p}}\sum_{(x,|y|)\in \Pi^t} \chi(a)\chi(x)\overline{\chi(y)},
$$
where the first summation is over all characters to modulus $p.$
Consider the sum $S = \sum_{a\in R}S(a)$ then
$$
S = \frac{\#R}{p-1}\sum_{\chi; U}\sum_{(x,|y|)\in \Pi^t}\chi(v)\chi(x)\overline{\chi(y)},
$$
where $\chi; U$ denotes the summation over all characters to prime modulus $p$  trivial on $U.$
It's now clear that
$$
|S| \leq \#R\frac{4p\log{p}}{(p-1)t} + 2\max_{\chi \neq 1}\left|
\sum_{(x,y)\in \Pi^t} \chi(x)\overline{\chi(y)} \right|.
$$
Using Lemma 1 we obtain the following estimate
$$
|S| \leq \frac{\#R}{4}
+ 20000 p^\frac{7}{8} \log^2{p} / \sqrt{t},
$$
therefore
$$
\frac{|S|}{\#R} < \frac 1 2,
$$
and the theorem follows.

\section{Proof of Theorem 2}
Define
$$
B(c) = \{(a,x)\bigm|\|\frac{ax}{p}\|\leq \frac{1}{cx}, a\in R, 1\leq x <p \}
$$
Let $B(c,c') = B(c)\diagdown B(c').$ Also we define a function $f_a(x)$ by the condition
$$
f_a(x) = 
\begin{cases}
c & \text{ if $(a,x)\in B(c,c+1)$ for some $c \in \mathbbm{N} $}\\
0 & \text{otherwise.}
\end{cases}
$$

Consider the sum
$$
S_a = \sum_{x=1}^{p-1} f_a(x).
$$

We distinguish two cases.

If $f_a(x) = c > 1,$ then 
$$
\frac{1}{(c+1)x} \leq \|\frac{ax}{p}\| < \frac{1}{cx}.
$$
There exists an integer $b$ such that
$$
\frac{1}{(c+1)x} \leq \left| \frac{ax}{p} - b \right| < \frac{1}{cx}.
$$
Hence
$$
\left| \frac{a}{p} - \frac{b}{x} \right| < \frac{1}{cx^2}
$$
and by Lemma A $\frac{b}{x}$ is convergent fraction to $\frac{a}{p}$.
Therefore $x=x_n$ is a denominator of a certain (say $n$-th) convergent (n=0, 1, $\cdots$).
By Lemma B
$$
\frac{1}{(a_{n+1}+2)x} \leq \|\frac{ax}{p}\| < \frac{1}{a_{n+1}x},
$$
therefore either $a_{n+1} = c,$ or $a_{n+1} = c - 1.$

If $f_a(x) = c = 1,$ then either $x$ is convergent's denominator and $a_{n+1} = 1,$ or $x$ is not a denominator to a convergent to $\frac{a}{p}.$

So we see that
$$
S_a = \sum b_i + \sum \delta_i + |W|,
$$
where $\delta_i \in \{-1,0\},$ $W\subset \left\{x \bigm| \left\| \frac{ax}{p} \right\| < \frac{1}{x} \right\}$

Therefore
$$
S_a \geq \sum_{i=1}^{l} b_i - 5\log{p},
$$
where $5\log{p}$ is an upper bound for the continued fraction's length.

Let $\Omega$ be the subset in $R$ such that all partial quotients of the elements of $\Omega$ are less   than $t = 16\log{p}.$
Hence, by Theorem 1, $\#\Omega > \#R/2.$
So if $a\in \Omega,$ then $f_a(x) < t.$
Hence by the partial summation
$$
\sum_{a\in \Omega } S_a \leq \sum_{c\leq t}c\cdot \#B(c,c+1) \leq
\sum_{c\leq t}\#B(c).
$$

Let's estimate $\#B(c)$.

It's clear that
\begin{equation}
\begin{array}{l}
\#B(c) \leq 2\cdot \# \{(b,x)\mid b<\frac{p}{cx}, b\in x\cdot K\}\\
\leq 2\cdot \#\{(b,x)\in \Pi^c\mid b\in x\cdot K\}\\
= 2\frac{\#R}{p-1} \sum_{\chi ; R}\sum_{(x,u)\in \Pi^c}\chi(v)\chi(u)\overline{\chi(x)},
\end{array}
\end{equation}
where the $\sum_{\chi;R}$ denotes the summation over characters $\chi$  trivial on $R$.

Note
$\#R|(p-1)$ and there exist exactly $ (p-1)/\# R$ trivial on $R$ characters.
Thus
$$
\#B(c) \leq 2\frac{\#R}{p-1}\# \Pi^c
+4\max_{\chi }\left| \sum_{(x,u)\in \Pi^c}
\chi(u)\overline{\chi(x)}\right|
$$
where maximum is taken over all non-principal characters to modulus $p$.
We can now use Lemma 1 to obtain an estimate
$$
\#B(c) \leq 4\frac{\#R}{p-1}\frac{p}{c}\log{p} + 40000 p^\frac{7}{8} \log^2{p} / \sqrt{c}.
$$
Therefore
$$
\sum_{a\in \Omega } S_a \leq 
190\cdot\#R\log{p}\log\log{p}
+ 8 \cdot 10^6 p^{7/8}\log^{5/2}{p} .
$$
Dividing by $\#\Omega>\#R/2$
we get
$$
\frac{1}{\#\Omega}\sum_{a\in \Omega } S_a \leq 400\log{p}\log\log{p}
$$
because of $\#R \geq 10^8 p^{7/8}\log^{5/2}{p}.$ Hence there exists an element $a$ in $\Omega$ such that $S_a \leq 400\log{p}\log\log{p}.$
Therefore there exists element $a$ in $\Omega$ such that
$$
\sum{b_i(a)} \leq 500 \log{p}\log\log{p}.
$$
Theorem 2 is proved.

\section{ Acknowledgements }
The authors are grateful to Prof. S.V. Konyagin for pointing out
an opportunity of improvement of Lemma 1 and hence the main result of the paper.

\end{document}